\documentclass[11pt,a4paper]{article}

\usepackage{amsmath,amsfonts,amssymb,latexsym,graphics,epsfig,url}
\usepackage{xcolor}
\usepackage{amsthm}
\usepackage[english]{babel}
\usepackage{mathdots}
\usepackage{graphicx}
\usepackage{soul} %subrrayar
\usepackage[utf8]{inputenc}
\usepackage{diagbox}
\usepackage[makeroom]{cancel}
\usepackage{float}              
\newtheorem{theorem}{Theorem}[section]
\newtheorem{proposition}[theorem]{Proposition}
\newtheorem{lemma}[theorem]{Lemma}
\newtheorem{corollary}[theorem]{Corollary}

\newtheorem{example}[theorem]{Example}

\DeclareMathOperator{\dist}{dist}

\DeclareMathOperator{\tr}{tr}
\DeclareMathOperator{\spec}{sp}
\DeclareMathOperator{\spann}{span}

\def\R{\ns R}
\def\b{\mbox{\boldmath $b$}}

\def\u{\mbox{\boldmath $u$}}

\def\vecv{\mbox{\boldmath $v$}}

\def\vec0{\mbox{\boldmath $0$}}

\def\A{\mbox{\boldmath $A$}}
\def\B{\mbox{\boldmath $B$}}

\def\E{\mbox{\boldmath $E$}}
\def\D{\mbox{\boldmath $D$}}

\def\I{\mbox{\boldmath $I$}}
\def\J{\mbox{\boldmath $J$}}

\def\L{\mbox{\boldmath $L$}}
\def\M{\mbox{\boldmath $M$}}

\def\X{\mbox{\boldmath $X$}}
\def\Y{\mbox{\boldmath $Y$}}

\def\I{\mbox{\boldmath $I$}}
\def\J{\mbox{\boldmath $J$}}

\def\R{\mbox{\boldmath $R$}}

\def\1{\mbox{\boldmath $1$}}

\def\Re{\mathbb R}

\begin{document}
	
	\title{On two algebras of token graphs
		%any graph (and of a cycle)
		\thanks{This research has been supported by
			AGAUR from the Catalan Government under project 2021SGR00434 and MICINN from the Spanish Government under project PID2020-115442RB-I00.
			The research of M. A. Fiol was also supported by a grant from the  Universitat Polit\`ecnica de Catalunya with references AGRUPS-2022 and AGRUPS-2023.}
	}
	\author{M. A. Reyes$^a$, C. Dalf\'o$^a$, M. A. Fiol$^b$\\
		\\
		{\small $^a$Dept. de Matem\`atica, Universitat de Lleida, Lleida/Igualada, Catalonia}\\
		{\small {\tt \{monicaandrea.reyes,cristina.dalfo\}@udl.cat}}\\
		%,visitant.arnau.messegue
		{\small $^{b}$Dept. de Matem\`atiques, Universitat Polit\`ecnica de Catalunya, Barcelona, Catalonia} \\
		{\small Barcelona Graduate School of Mathematics} \\
		{\small  Institut de Matem\`atiques de la UPC-BarcelonaTech (IMTech)}\\
		{\small {\tt miguel.angel.fiol@upc.edu} }\\
	}

	\date{}
	\maketitle
	
	\begin{abstract}
		The $k$-token graph $F_k(G)$ of a graph $G$ is the graph whose vertices are the $k$-subsets of vertices from $G$, two of which being adjacent whenever their symmetric difference is a pair of adjacent vertices in $G$.
		In this article, we describe some properties of the Laplacian matrix $\L_k$ of $F_k(G)$ and the Laplacian matrix $\overline{\L}_k$ of the $k$-token graph $F_k(\overline{G})$ of its complement $\overline{G}$.
		In this context, a result about the commutativity of the matrices $\L_k$ and $\overline{\L}_k$ was given in [C. Dalf\'o, F. Duque, R. Fabila-Monroy, M. A. Fiol, C. Huemer, A. L. Trujillo-Negrete, and F. J. Zaragoza Mart\'{\i}nez, 
		On the Laplacian spectra of token graphs, 
		{\em Linear Algebra Appl.} {\bf 625} (2021) 322--348], but the proof was incomplete, and there were some typos. Here we give the correct proof.
		Based on this result, and fixed the pair $(n,k)$ and the graph $G$, we first introduce a `local' algebra ${\cal L}(G)$, generated by the pair $(\L_k, \overline{\L}_k)$, showing its closed relationship with the Bose-Mesner algebra of the Johnson graphs $J(n,k)$.
		Finally, fixed only $(n,k)$, we present a `global' algebra ${\cal A}(n,k)$ that contains ${\cal L}(G)$ together with the Laplacian and adjacency matrices of the $k$-token graph of any graph $G$ on $n$ vertices.
	\end{abstract}
	
	
	
	%%%%%%%%%%%%%%%%%%%%%%%%%%%%%%%%%%%%%%%%%%%%%%%%%%%%%%%%%%%%%%%%%%%%%%%%%%%%%%%%%%%%%%%%%%%%%%%%%%%%%%%%%%%%%%%%%%%%%%%%%%%%%%%%%%%%
	
	\section{Introduction}
	
	In graph theory, there are different operations that construct a `large' graph from a `smaller' one. Then, the question is what properties of the former can be deduced (or, at least, approximate) from the properties of the latter. One of these operations that recently received some attention in the literature is the construction of token graphs. Namely, from a given graph $G$ on $n$ vertices and an integer $k<n$, we construct its token graph $F_k(G)$ on ${n\choose k}$ vertices. 
	The name `token graph' comes from an observation in
	Fabila-Monroy, Flores-Pe\~{n}aloza,  Huemer,  Hurtado,  Urrutia, and  Wood \cite{ffhhuw12}, that vertices of $F_k(G)$ correspond to configurations of $k$ indistinguishable tokens placed at distinct vertices of $G$, where two configurations are adjacent whenever one configuration can be reached
	from the other by moving one token along an edge from its current position to an unoccupied vertex. For example, Figure \ref{fig2} shows the 2-token of the complete graph $K_5$.
	The $k$-token graphs are also called symmetric $k$-th power of graphs in Audenaert, Godsil, Royle, and Rudolph \cite{agrr07}, and $k$-tuple vertex graphs in Alavi, Lick, and Liu \cite{all02}.
	Token graphs have some applications in physics. For instance, a relationship between token graphs and the exchange of Hamiltonian operators in
	quantum mechanics is given in Audenaert, Godsil, Royle, and Rudolph \cite{agrr07}. Our interest in the study of token graphs is motivated by some of their applications in %both physics and 
	mathematics and computer science: Analysis of complex networks, coding theory, combinatorial designs (by means of Johnson graphs),
	%quantum mechanics, %(ref. [2] of the paper),
	algebraic graph theory, enumerative combinatorics, the study of symmetric functions, etc.

	We describe some interesting properties of the Laplacian matrices of $F_k(G)$ and $F_k(\overline{G})$, where $\overline{G}$ is the complement graph of $G$. This leads us to reveal the closed relationship between the algebra of the pair $(F_k(G), F_k(\overline{G}$)) with the Bose-Mesner algebra of the Johnson graphs $J(n,k)$.
	This paper is structured as follows. In Section \ref{resultats-coneguts}, we give the preliminaries and background together with some known results. In particular, we introduce the family of Laplacian predistance polynomials, deriving a simple Hoffman-like condition for connectedness. The important family of distance-regular graphs known as Johnson graphs are also introduced as token graphs, together with some of their properties.  In Section \ref{L-spec-token}, we first recall some known properties of token graphs' Laplacian spectra and obtain new related results. Moreover, a result concerning the commutativity of the Laplacian matrices of a graph and its complement was given in \cite{ddffhtz21}, but the proof was incomplete, and there were some typos. Here we give the correct proof.
	Sections \ref{local-algebra-token} and \ref{global-algebra-token} contain the main results of this paper by introducing two new algebras for token graphs, which are closely related to the Bose-Mesner algebra of Johnson graphs. More precisely, in Section \ref{local-algebra-token}, for fixed values  $(n,k)$ and graph $G$, the `local' algebra ${\cal L}(G)$ is generated by the pair $(\L_k, \overline{\L}_k)$, and contains  the Bose-Mesner algebra ${\cal A}(J(n,k))$.
	Finally, fixed only $(n,k)$, the `global' algebra ${\cal A}(n,k)$ is generated by some `elementary' adjacency matrices, and contains ${\cal L}(G)$ together with the Laplacian and adjacency matrices of the $k$-token graph of any graph $G$ on $n$ vertices.

	%In Figures \ref{F_2(C_9)}, \ref{F_2(C_{10})}, and  \ref{F_2(C_8)}, we show the 2-token graphs of cycles $C_9$, $C_{10}$,  and $C_8$, respectively
	%%%

	%%%%%%%%%%%%%%%%%%%%%%%%%%%%%%%%%%%%%%%%%%%%%%%%%%%%%%%%%%
	%%%%%%%%%%%%%%%%%%%%%%%%%%%%%%%%%%%%%%%%%%%%%%%%%%%%%%%%%
	
	\section{Preliminaries and background}
	\label{resultats-coneguts}
	
	%\subsection{Some notation and basic facts}
	
	%%%%%%%%%%%%%%%%%%%%%%%%%%%%%%%%%%%%%%%%%%%%%%%%%%%%%%%%%%%%%%%%%%
	
	In this section, we give the basic notation and definitions. 
	Moreover, together with some known results, we also give some basic new results. Let us start with the formal definition of token graphs.
	
	\subsection{Token graphs}
	Let $G=(V,E)$ be a simple graph with vertex set $V=V(G)=\{1,2,\ldots,n\}$ and edge set $E=E(G)$. By convenience, we consider every edge $e=\{u,v\}$ constituted by two opposite arcs $(u,v)$ and $(v,u)$. Let $N(u)$ denote the set of vertices adjacent to $u\in V$, so that the minimum degree of $G$ is $\delta(G)=\min_{u\in V}|N(u)|$. For a given integer $k$ such that $1\leq k \leq n$, the \textit{$k$-token graph} $F_k(G)$ of $G$ is the graph whose vertex set $V (F_k(G))$ consists of the ${n \choose k}$ $k$-subsets of vertices of $G$, and two vertices $A$ and $B$ of $F_k(G)$ are adjacent if and only if their symmetric difference $A \bigtriangleup B$ is a pair $\{a,b\}$ such that $a\in A$, $b\in B$, and $(a,b)\in E(G)$.
	Then, if $G$ has $n$ vertices and $m$ edges, $F_k(G)$ has ${n \choose k}$ vertices and ${n-2 \choose k-1}m$ edges. (Indeed, for each edge of $G$, there are 
	${n-2 \choose k-1}$ edges in $F_k(G)$.)
	We also use the notation $\{a,b\}$, with $a, b\in V$, for a vertex of a 2-token graph. Moreover, we use $ab$ for the same vertex in the figures. 
	For convenience, we define $F_0(G)$ as a singleton
	$\{u\}$.
	Moreover, $F_1(G)=G$ and, by symmetry,
	$F_k(G)=F_{n-k}(G)$.
	Moreover, if $G$ is bipartite, so it is $F_k(G)$ for any $k=1,\ldots,|V|-1$.
	
	\subsection{Graph spectra and orthogonal polynomials}
	
	Given a graph $G$
	with adjacency matrix $\A$ and spectrum
	$$
	\spec \A=\{\theta_0^{m_0},\theta_1^{m_1},\ldots,\theta_d^{m_d}\},
	$$
	where $\theta_0>\theta_1>\cdots>\theta_d$,
	the {\em predistance polynomials} $p_0,p_1,\ldots,p_d$, introduced by Fiol and Garriga in \cite{fg97}, are a sequence of orthogonal polynomials, such that $\deg p_i=i$, with respect to the scalar product
	\begin{equation}
		\langle f,g\rangle_{A}=\frac{1}{n}\tr (f(A)g(A))=
		\frac{1}{n}\sum_{i=0}^d m_i f(\theta_i)g(\theta_i),
		\label{esc-prod}
	\end{equation}
	normalized in such a way that
	$$
	\|p_i\|_{A}^2=p_i(0).
	$$
	Notice that, using this norm, $\|1\|_A=1$.
	We can apply the Gram-Schmidt process to $1,x,x^2,\ldots,x^d$ to find these polynomials and normalize accordingly.
	The name given to them is justified because, when $G$ is a distance-regular graph, the predistance polynomials become the well-known {\em distance polynomials} that applied to $\A$ yield the distance matrices, that is, $p_i(\A)=\A_i$.
	%%%%%%%%%%%%%%%%%%%%%%%%%%%%%%%%%%%%%%%%%%%%%%%%
	
	In a similar way, we can introduce another sequence of orthogonal polynomials from the spectrum of the Laplacian matrix.
	Let $G$ be a graph on $n$ vertices, with Laplacian matrix $\L=\D-\A$ and spectrum
	$$
	\spec \L= \{\lambda_0^{m_0},\lambda_1^{m_1},\ldots, \lambda_l^{m_l}\},
	$$
	where $\lambda_0(=0)<\lambda_1<\cdots <\lambda_l$.
	The Laplacian predistance polynomials $q_0,q_1,\ldots,q_l$ are a sequence of orthogonal polynomials, such that $\deg q_i=i$, with respect to the scalar product
	\begin{equation}
		\langle f,g\rangle_L=\frac{1}{n}\sum_{i=0}^l m_if(\lambda_i)g(\lambda_i),
		\label{L-esc-prod}
	\end{equation}
	normalized in such a way that
	$$
	\|q_i\|^2_L=q_i(0).
	$$
	As in the previous case, the norm associated with \eqref{L-esc-prod} satisfies $\|1\|_L=1$.
	%%%%%%%%%%%%%%%%%%%%%%%%%%%%%%%%%%
	In the following proposition, we use the (adjacency and Laplacian) predistance polynomials to give two basic characteristics of graphs, namely, regularity and connectivity.
	The first assertion is a reformulation of the celebrated Hoffman characterization of regularity, see \cite{h63}.
	The second statement is a Hoffman-type result obtained by using the Laplacian matrix (see, for example, Arsi\'{c}, Cvetkovi\'{c}, Simi\'{c}, and \v{S}kari\'{c} \cite[p. 19]{acss12} or 
	Van Dam and Fiol \cite[p. 247]{vf14}).
	\begin{proposition}
		Let $G$ be a graph with adjacency matrix $\A$ with $d+1$ distinct eigenvalues and Laplacian matrix $\L$ with $l+1$ distinct eigenvalues. Let $p_0,\ldots,p_d$ be its predistance polynomials, and $q_0,\ldots,q_l$ its  Laplacian predistance polynomials. Consider the sum polynomials
		$H=\sum_{i=0}^d p_i$ and  $H_L=\sum_{i=0}^l q_i$. Let $\J$ the all-$1$ matrix. Then, the following statements hold.
		\begin{itemize}
			\item[$(i)$] \cite{h63}
			$G$ is connected and regular if and only if $H(\A)=\J$. 
			\item[$(ii)$] \cite{acss12,vf14}
			$G$ is connected if and only if $H_L(\L)=\J$.
		\end{itemize}
	\end{proposition}

	%%%%%%%%%%%%%%%%%%%%%%%%%%%%%%%%%%%%%%%%%%%%%%%%%%%%%%%%%%%%%%%%%%%%%%%
	
	\subsection{Distance-regular graphs and the Bose-Mesner algebra}
	
	Recall that a graph $G$ with diameter $D=d$ is \textit{distance-regular} if, for any pair of vertices $u,v$, the intersection parameters $a_i(u)=|G_i(u)\cap G(v)|$, $b_i(u)=|G_{i+1}(u)\cap G(v)|$, and $c_i(u)=|G_{i+1}(u)\cap G(v)|$, for $i=0,1,\ldots,d$, only depend on the distance $\dist(u,v)=i$, where $G(v)$ is the set of vertices adjacent to vertex $v$ and $G_i(u)$ is the set of vertices at distance $i$ from $u$.
	Then, we have the intersection array
	\begin{equation}
		\label{intersec-array}
		\iota(G)=\left\{
		\begin{array}{ccccc}
			- & c_1 & \cdots & c_{d-1} & c_{d}\\
			a_0 & a_1 & \cdots & a_{d-1} & a_{d}\\
			b_0 & b_1 & \cdots & b_{d-1} & -
		\end{array} \right\}
	\end{equation}
	or quotient matrix
	\begin{align}
		\B & =\left(
		\begin{array}{ccccc}
			a_0 & c_1 &        &        &      \\
			b_0 & a_1 & c_2    &        &      \\
			& b_1 & a_2    & \ddots &      \\
			&     & \ddots & \ddots &  c_{d} \\
			&     &        & b_{d-1}&  a_{d} \\
		\end{array}
		\right).\label{A-quotient}
	\end{align}
	
	In terms of the predistance polynomials, a characterization of distance-regularity is the following.
	A graph $G$ with diameter $d$ and distance matrices
	$$
	\A_0(=I), \A_1(=\A),\ldots, \A_{d}
	$$
	is \textit{distance-regular} if and only if the predistance polynomials
	$p_0,p_1,\ldots,p_{d}$ satisfy
	\begin{equation}
		\A_i=p_i(\A),\qquad i=0,1,\ldots,{d}.
		\label{pi(A)=Ai}
	\end{equation}
	In fact, as it was proved in Fiol, Garriga, and Yebra \cite{fgy96}, the last equality in \eqref{pi(A)=Ai} involving the highest degree polynomial suffices, namely, $\A_d=p_d(\A)$.
	
	%%%%%%%%%%%%%%%%%%%%%%%%%%%%%%%%%%%%%%%%%%%%%%%%%%%%%%%%%%
	
	%\subsection*{The Bose-Mesner algebra}
	A way of introducing the Bose-Mesner algebra of a graph is as follows.
	Let $G$ be a graph with diameter $D$, distance matrices $\A_0(=\I), \A_1(=\A), \A_2,\ldots, \A_D$, and $d+1$ distinct eigenvalues.
	Consider the vector spaces
	\begin{align*}
		{\cal A} & = \spann\{\I,\A,\A^2,\ldots,\A^d\},\\
		{\cal D} &= \spann\{\I, \A, \A_2,\ldots,\A_D\},
	\end{align*}
	with dimensions $d+1$ and $D+1$, respectively.
	Then,
	${\cal A}$  is an algebra with the ordinary product of matrices,
	known as the \textit{adjacency algebra} of $G$.
	Moreover,
	${\cal D}$ is an algebra with the entry-wise (or Hadamard
	product) of matrices, defined by $(\X \circ \Y)_{uv}= \X_{uv}\Y_{uv}$, called the \textit{distance $\circ$-algebra} of $G$.
	Notice that, if $G$ is regular, then 
	$\I,\A,\J\in {\cal A} \cap {\cal D}$.
	The most interesting case happens when both algebras coincide.
	\begin{theorem}
		Let $G$, ${\cal A}$, and ${\cal D}$ be as above. Then, $G$ is distance-regular if and only if
		\begin{equation}
			\label{A=D}
			{\cal A}={\cal D},
		\end{equation}
		so that ${\cal A}$ is an algebra with both the ordinary product and the Hadamard product of matrices, and it is known as the {\em Bose-Mesner algebra}  associated with $G$,
	\end{theorem} 
	Since $D\le d$, the equality in \eqref{A=D}  is equivalent to $\dim({\cal A} \cap {\cal D})=d+1.$
	
	\subsection{Johnson graphs}
	
\begin{figure}[t]
		\begin{center}
			\includegraphics[width=5cm]{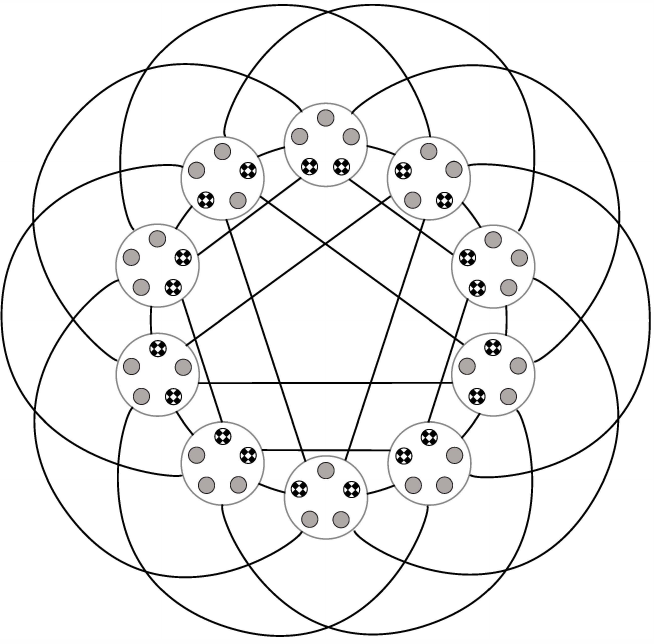}
			%			$\mathrm{Sp}\, G=\{0,2,3,4,5\}\qquad \mathrm{Sp}\,F_2(G)=\{0,2,3^2,4,5^3,7,8\}$.
			\caption{The Johnson graph $J(5,2)=F_2(K_5)=\overline{P}$.}
			\label{fig2}
		\end{center}
	\end{figure}

	An important class of token graphs are the {\em Johnson graphs} $J(n,k)$ (see Fabila-Monroy, Flores-Pe\~{n}aloza, Huemer, Hurtado, Urrutia, and 
	Wood~\cite{ffhhuw12}), that correspond to the $k$-token of the complete graph
	$F_k(K_n)$. In particular $J(n,1)=K_n$, the graph $J(4,2)$ is the octahedron, and $J(5,2)$ is the complement of the Petersen graph $P$, see Figure \ref{fig2}.
	
	In general, the Johnson graph $J(n,k)$
	is distance-transitive (and, hence, distance-regular) satisfying the following properties:
	\begin{itemize}
		\item[{\bf P1.}]
		The Johnson graph $J(n,k)$, with $k\le n-k$, is a \textit{distance-regular} graph with degree $k(n-k)$, diameter $d=k$, and intersection parameters
		$$
		b_j=(k-j)(n-k-j),\quad c_j=j^2,\quad  \mbox{for }j=0,1,\ldots,d.
		$$
		\item[{\bf P2.}]
		The Laplacian eigenvalues $\lambda_j$ (and their multiplicities $m_j$) of $J(n,k)$ are
		$$
		\lambda_j=j(n+1-j),\quad m_j={n\choose j}-{n\choose j-1},\quad \mbox{for }j=0,1,\ldots,k.
		$$
		For instance, the Laplacian eigenvalues of $J(n,4)$ are
		$0, n, 2(n-1), 3(n-2)$ and $4(n-3)$.
		\item[{\bf P3.}]
		Any pair of vertices are at distance $j$, with $0\le j\le k$, if and only if they share $k-j$ elements in common.
		\item[{\bf P4.}]
		The Johnson graph $J(n,k)$ is maximally connected, that is, $\kappa=k(n-k)$.
	\end{itemize}
	
	%%%%%%%%%%%%%%%%%%%%%%%%%%%%%%%%%%%%%%%%%%%%%%%%%%%%%%%
	\section{On the Laplacian spectra of token graphs}
	\label{L-spec-token}
	
	This section deals with some properties of the Laplacian spectra of token graphs.
	With this aim, let $G$ be a graph on $n$ vertices. Let $[n]:=\{1,\ldots,n\}$ and ${[n]\choose k}$ denote the set of $k$-subsets of $[n]$, which is the set of vertices of the $k$-token graph $F_k(G)$.
	For our purpose, it is convenient to denote by $W_n$ the set of all column vectors  $\vecv$ with $n$ entries such that $\vecv^{\top }\1 = 0$.
	Let $\lambda_1,\lambda_2, \ldots, \lambda_n$  be the eigenvalues of the Laplacian matrix $\L(G)$ of a graph $G$, with $(0=)\lambda_1\le \lambda_2\le \cdots \le \lambda_n$. The second smallest eigenvalue $\lambda_2$ is known as the \textit{algebraic connectivity} $\alpha(G)$.
	
	%%%%%%%%%%%%%%%%%%%%%%%%%%%%%%%%%%%%%%%%%%%%%%%%%%%%%%%%%%
	%\textbf{The $(n,k)$-binomial matrix}
	
	Given integers $n$ and $k(\le n)$, the $(n;k)$-\textit{binomial matrix} $\B$ is an ${n \choose k}\times n$ matrix whose rows are the characteristic vectors of the $k$-subsets of $[n]=\{1,\ldots,n\}$ in a given order. Thus, if the $i$-th $k$-subset is $A$, then
	$$
	(\B)_{ij}=
	\left\lbrace
	\begin{array}{ll}
		1 & \mbox{if } j\in A,\\
		0 & \mbox{otherwise.}
	\end{array}
	\right.
	$$
	For instance, for $n=4$ and $k=2$, we have
	$$
	\B=\left(
	\begin{array}{cccc}
		1 & 1 & 0 & 0\\
		1 & 0 & 1 & 0\\
		1 & 0 & 0 & 1\\
		0 & 1 & 1 & 0\\
		0 & 1 & 0 & 1\\
		0 & 0 & 1 & 1
	\end{array}
	\right).
	$$
	By a simple counting argument, one can check that this matrix satisfies
	$$
	\B^{\top}\B={n-2\choose k-1}\I + {n-2\choose k-2}\J,
	$$
	where $\J$ is the all-1 matrix, see Dalf\'o, Duque, Fabila-Monroy, Fiol, Huemer, Trujillo-Negrete, and Zaragoza Mart\'inez~\cite{ddffhtz21}.
	
	The following result enumerate other important properties of the $(n,k)$-binomial matrix, proved in the same paper.
	
	%%%%%%%%%%%%%%%%%%%%%%%%%%%%%%%%%%%%%%%%%%%%%%%%%%%%%%%%%%
	
	%\textbf{The Laplacian spectra of token graphs}
	\begin{theorem}[\cite{ddffhtz21}]
		\label{coro:LkL1}
		Let $G$ be a graph with Laplacian matrix $\L_1$. Let $F_k=F_k(G)$ be its token graph with Laplacian $\L_k$. Then, the following statements hold:
		\begin{itemize}
			\item[$(i)$]
			$\B\L_1=\L_k\B$.
			\item[$(ii)$]
			$\L_1=(\B^{\top}\B)^{-1}\B^{\top}\L_k\B=\frac{1}{{n-2\choose k-1}}\B^{\top}\L_k\B$.
			\item[$(iii)$]
			The column space (and its orthogonal complement) of $\B$ is $\L_k$-invariant.
			\item[$(iv)$]
			The characteristic polynomial of $\L_1$ divides the characteristic polynomial of $\L_k$.
			Thus, $\spec \L_1 \subseteq \spec \L_k$.
			\item[$(v)$]
			If $\vecv$ is a $\lambda$-eigenvector of $\L_1$, then $\B\vecv$ is a $\lambda$-eigenvector of $\L_k$.
			\item[$(vi)$]
			If $\u$ is a $\lambda$-eigenvector of $\L_k$ such that $\B^{\top}\u\neq \vec0$, then $\B^{\top}\u$
			is a $\lambda$-eigenvector of $\L_1$.
		\end{itemize}
	\end{theorem}
	
	%%%%%%%%%%%%%%%%%%%%%%%%%%%%%%%%%%%%%%%%%%%%%%%%%%%%%%%%%%%%%%%%%%%%%%%%
	
	%\blue{
		Concerning the Laplacian matrices of a graph and its complement, the following result was given in \cite{ddffhtz21}, but the proof was incomplete, and there were some typos. So, for completeness, we give here the correct proof.
		
		\begin{proposition}[\cite{ddffhtz21}]
			\label{propo:commute}
			Let $G=(V,E)$ be a graph on $n=|V|$ vertices, and let $\overline{G}$ be its complement. For a given $k$, the Laplacian matrices of their $k$-token graphs $\L_k=\L(F_k(G))$ and $\overline{\L}_k= \L(F_k(\overline{G}))$ commute:
			$$
			\L_k\overline{\L}_k=\overline{\L}_k\L_k.
			$$
		\end{proposition}
		
		\begin{proof}
			Let $\L=\L_k$ and $\overline{\L}=\overline{\L}_k$. We want to prove that $(\L \overline{\L})_{AB}=(\overline{\L}\L)_{AB}$ for every pair of vertices $A,B\in {[n]\choose k}$ of $F_k(G)$ and $F_k(\overline{G})$, respectively. To this end, we consider the different possible values of $|A\cap B|$.
			First, note that $(\L \overline{\L})_{AB}=(\overline{\L}\L)_{AB}=0$ when $|A\cap B|<k-2$ (see after \eqref{LnoL(AB)} the interpretation of the entry  $(\L \overline{\L})_{AB}$). Thus, we only need to consider the following three cases:
			\begin{enumerate} 
				\item[{\bf 1.}] 
				If $|A\cap B|=k$, that is $A=B$, we have 
				\begin{equation}
					\label{LLC(AB)1}
					(\L \overline{\L})_{AA}=(\overline{\L} \L)_{AA}=\deg_{F_k(G)}(A)\cdot\deg_{F_k(\overline{G})}(A).
				\end{equation}
				\item[{\bf 2.}]  
				If $|A\cap B|=k-1$, we can assume that $A=A'\cup \{a\}$ and $B=A'\cup \{b\}$, where $|A'|=k-1$ and $a,b\in[n]$ for $a\neq b$.
				Moreover, without loss of generality, we can assume that $a\sim b$ in $G$ implies that $a\not\sim b$ in $\overline{G}$. (If not, interchange the roles of $A$ and $B$.)
				Then, the different terms of the sum
				\begin{equation}
					\label{LnoL(AB)}
					(\L \overline{\L})_{AB}= \sum_{X\in {[n]\choose k}} (\L)_{AX}(\overline{\L})_{XB}
				\end{equation}
				can be seen as `walks' of length two, where the first step is done in $F_k(G)$, and the second step is done in $F_k(\overline{G})$. When we only have one step, the other corresponds to a loop in the initial or final vertex, represented as $A \rightarrow A$ or $B \rightarrow B$, respectively. Each step gives $-1$, except for each loop that provides the vertex degree. We multiply the giving of both steps. This yields the following contributions to the sum:
				\begin{enumerate}
					\item[$(i)$]
					First `add $b$', and then `delete $a$':
					$$
					A=A'\cup \{a\} \quad \rightarrow \quad X=(A\setminus \{c\})\cup\{b\} \quad \rightarrow \quad (X\setminus \{a\}) \cup \{c\}=B,
					$$
					where   $c\in A'$, $c\sim b$ in $G$; and $c\sim a$ in $\overline{G}$. Thus, for each $c\in N_G(b)\cap N_{\overline{G}}(a)\cap A'$, we get a term $(-1)(-1)=1$.
					\item[$(ii)$]
					First `delete $a$', and then  `add $b$':
					$$
					A=A'\cup \{a\} \quad \rightarrow \quad X=A'\cup\{c\}% , \ h\in
					\quad \rightarrow \quad A' \cup \{b\}=B,
					$$
					where $c\in \overline{A\cup B}$, $c\sim a$ in $G$; and $c\sim b$ in $\overline{G}$. Thus, for each $c\in N_G(a)\cap N_{\overline{G}}(b)\cap \overline{A\cup B}$, we get a term $(-1)(-1)=1$.
					\item[$(iii)$]
					First `change $a$ by $b$', and then  `keep $b$':
					$$
					A=A'\cup \{a\} \quad \rightarrow \quad X=A'\cup\{b\}=B \quad \rightarrow \quad B,
					$$
					where $a\sim b$ in $G$. This corresponds to the case when $X=B$ and $A\sim B$ in $F_k(G)$. Then, since the second step corresponds to the diagonal entry of $(\overline{\L})_{BB}$, this gives the term $(-1)\cdot\deg_{F_k(\overline{G})} B$.
					\item[$(iv)$]
					The other way around (first  `keep $a$', and then `change $a$ by $b$') corresponds to the case when $X=A$, but $A\not\sim B$ in $F_k(\overline{G})$ (since $a\not\sim b$ in $\overline{G}$). Then, this gives zero.
				\end{enumerate}
				To summarize, since all the other possibilities give a zero term, we conclude that, when $|A\cap B|=k-1$, the total value of the sum in \eqref{LnoL(AB)} is
				\begin{equation}
					\label{LLC(AB)2}
					(\L\overline{\L})_{AB}=|N_G(b)\cap N_{\overline{G}}(a)\cap A'|+|N_G(a)\cap N_{\overline{G}}(b)\cap \overline{A\cup B}|-\deg_{F_k(\overline{G})} (B).
				\end{equation}
				With respect to  $(\overline{\L}\L)_{AB}$, note that, since the involved matrices are symmetric, we can compute
				$$
				(\overline{\L}\L)_{BA}=\sum_{X\in{[n]\choose k}}(\overline{\L})_{BX}(\L)_{XA}.
				$$
				Reasoning as before, the walks $(i)$ and $(ii)$ (from $A$ to $B$) become the walks from $B$ to $A$ if we interchange these two vertices, the vertices $a$ and $b$, and the graphs $G$ and $\overline{G}$.
				Then, the cases $(i)$--$(iv)$ become:
				\begin{itemize}
					\item[$(i')$]
					First `add $a$', and then  `delete $b$':
					$$
					B=A'\cup \{b\} \quad \rightarrow \quad X=(B\setminus \{c\})\cup\{a\} \quad \rightarrow \quad (X\setminus \{b\}) \cup \{c\}=A,
					$$
					where $c\in A'$, $c\sim a$ in $\overline{G}$, and $c\sim b$ in $G$. Thus, for each $c\in N_G(b) \cap N_{\overline{G}}(a) \cap A'$, we get again a term $(-1)(-1)=1$.
					\item[$(ii')$]
					First `delete $b$', and then  `add $a$':
					$$
					B=A'\cup \{b\} \quad \rightarrow \quad X=A'\cup\{c\}% , \ h\in
					\quad \rightarrow \quad A' \cup \{a\}=A,
					$$
					where $c\in \overline{A\cup B}$, $c\sim b$ in $\overline{G}$, and $c\sim a$ in $G$. Thus, for each $c\in N_G(a)\cap N_{\overline{G}}(b)\cap \overline{A\cup B}$, we get a term $(-1)(-1)=1$.
					\item[$(iii')$]
					First  `keep $b$', and then  `change $b$ by $a$':
					$$
					B=A'\cup \{a\} \quad \rightarrow \quad X=B \quad \rightarrow \quad (B\setminus \{b\})\cup\{a\}=A,
					$$
					where $a\sim b$ in $G$. The first step corresponds to the diagonal entry of  $(\overline{\L})_{BB}$, which gives the term $(-1)\cdot\deg_{F_k(\overline{G})} (B)$ again.
					\item[$(iv')$]
					The other case (first `change $b$ by $a$', and then `keep $a$') corresponds to $X=A$, and, since $A\not\sim B$, this  gives zero.
				\end{itemize}
				Consequently,  $(\overline{\L}\L)_{AB}=(\overline{\L}\L)_{BA}$ equals the same expression in \eqref{LLC(AB)2}.
				\item [{\bf 3.}]
				If $|A\cap B|=k-2$, we can assume that $A=A'\cup \{a,b\}$ and $B=A'\cup \{c,d\}$, where $|A'|=k-2$ and $a,b,c,d\in[n]$ with $\{a,b\}\cap\{c,d\}=\emptyset$.
				Now, the non-zero entries of $(\L\overline{\L})_{AB}$ must correspond to some of the following cases:
				\begin{itemize}
					\item [$(1)$]
					$a\sim c$ in  $G$  and $b\sim d$ in  $\overline{G}$.
					\item [$(2)$]
					$a\sim d$ in  $G$ and $b\sim c$ in  $\overline{G}$.
					% \end{itemize}
				% \begin{itemize}
					\item  [$(3)$]
					$b\sim d$ in  $G$  and $a\sim c$ in  $\overline{G}$.
					\item  [$(4)$]
					$b\sim c$ in  $G$    and $a\sim d$ in  $\overline{G}$.
				\end{itemize}	
				Note $(1)$ and $(3)$ exclude each other, and the same happens with $(2)$ and $(4)$. 
				Then, at most two of these four cases can be satisfied. 
				Let us assume that we are in some of the cases $(1)$ and $(2)$ (the other cases can be dealt with similarly). Therefore, we get the following situations:
				\begin{enumerate}
					\item[$(i)$]
					First $a\rightarrow c$, and then $b\rightarrow d$:
					$$
					A=A'\cup \{a,b\} \quad \rightarrow \quad X=A'\cup\{c,b\} \quad \rightarrow \quad A'\cup\{c,d\}=B.
					$$
					\item[$(ii)$]
					First $c\rightarrow b$, and then $d\rightarrow a$:
					$$
					B=A'\cup \{c,d\} \quad \rightarrow \quad X=A'\cup\{b,d\} \quad \rightarrow \quad A'\cup\{b,a\}=A.
					$$
					In each case, we have the term $(-1)(-1)=1$. (Notice that, possibly, both cases hold, and then $(\L\overline{\L})_{AB}=2$.)
				\end{enumerate}
				Similarly, the term  $(\overline{\L}\L)_{AB}=(\overline{\L}\L)_{BA}$ can be computed as:
				\begin{enumerate}
					\item[$(i')$]
					First $d\rightarrow b$, and then $c\rightarrow a$:
					$$
					B=A'\cup \{c,d\} \quad \rightarrow \quad X=A'\cup\{c,b\} \quad \rightarrow \quad A'\cup\{a,b\}=A.
					$$
					\item[$(ii')$]
					First $b\rightarrow d$, and then $a\rightarrow c$:
					$$
					B=A'\cup \{c,d\} \quad \rightarrow \quad X=A'\cup\{a,d\} \quad \rightarrow \quad A'\cup\{c,d\}=B.
					$$
					Thus, in each case, we have again the term $(-1)(-1)=1$.
				\end{enumerate}
			\end{enumerate}
			Finally, considering all the cases {\bf 1}--{\bf 3}, we have that $(\L\overline{\L})_{AB}=(\overline{\L}\L)_{AB}$ for any $A$ and $B$, as claimed.
		\end{proof}
		Here, it is worth mentioning that, in general, this result does NOT hold for the respective adjacency matrices of $F_k(G)$ and $F_k(\overline{G})$. 
		%}
	
	By a theorem of Frobenius, we get the following consequence.
	\begin{corollary}
		The ${n\choose k}$ eigenvalues of $\L_k$ and $\overline{\L}_k$ can be matched up as
		$\lambda_i\leftrightarrow \overline{\lambda}_i$ in such a way that the $n$ eigenvalues of any polynomial $p(\L_k , \overline{\L}_k )$ in the two matrices is the multiset of the values $p(\lambda_i,\overline{\lambda}_i)$.
	\end{corollary}
	Moreover, since the Laplacian matrix of $J(n,k)$ is $\L_J=\L_k+\overline{\L}_k$,
	$\L_J$ commutes with both $\L_k$ and $\overline{\L}_k$, and
	every eigenvalue  $\lambda_J$ of $J(n,k)$ is the sum of one eigenvalue $\lambda$ of $F_k(G)$ and one eigenvalue $\overline{\lambda}$ of $F_k(\overline{G})$. More precisely, we can state the following result about how to pair the eigenvalues.
	\begin{proposition}[\cite{df22}]
		\label{pro:pairing}
		Let $\L_k$ and $\overline{\L}_k$ be the Laplacian matrices of $F_k(G)$ and $F_k(\overline{G})$, respectively.
		For %$i=0,\ldots,k$ and
		$j=0,1,\ldots,k$,
		let $\lambda_j=j(n+1-j)$ and $m_j={n\choose j}-{n\choose j-1}$ be the eigenvalues and multiplicities of $J(n,k)$.
		Let $\lambda_{j1},\lambda_{j2},\ldots,\lambda_{jm_j}$ be
		the eigenvalues in $\spec F_j(G)\setminus \spec F_{j-1}(G)$, with $\lambda_{j1}\le \lambda_{j2}\le \cdots \le\lambda_{jm_j}$ 
		(that is, the {\em non-trivial} eigenvalues of $F_j(G)$).
		Let $\overline{\lambda}_{j1},\overline{\lambda}_{j2},\ldots,\overline{\lambda}_{jm_j}$ be
		the eigenvalues in $\spec F_j(\overline{G})\setminus \spec F_{j-1}(\overline{G})$, with $\overline{\lambda}_{j1}\ge \overline{\lambda}_{j2}\ge \cdots \ge\overline{\lambda}_{jm_j}$.
		% (that is the \blue{non-trivial} eigenvalues of $F_j(\overline{G})$).
		Then,
		\begin{equation}
			\lambda_{jr}+\overline{\lambda}_{jr}=\lambda_j\qquad \mbox{for } r=0,1,\ldots,m_j.
			\label{eq:pairing}
		\end{equation}
	\end{proposition}
	
	Part of this result was obtained in Dalf\'o, Duque, Fabila-Monroy, Fiol, Huemer, Trujillo-Negrete, and Zaragoza Mart\'{\i}nez \cite{ddffhtz21}, but they did not specify each pair of eigenvalues of $F_k(G)$ and $F_k(\overline{G})$ that yields the corresponding eigenvalue of $J(n, k)$. Later, Dalf\'o and Fiol addressed this question and established such a pairing %in Lemma 3.1 
	in \cite[Lemma 2.3]{df22}, see the following examples.
	
\begin{figure}[!ht]
		\begin{center}
			\includegraphics[width=6cm]{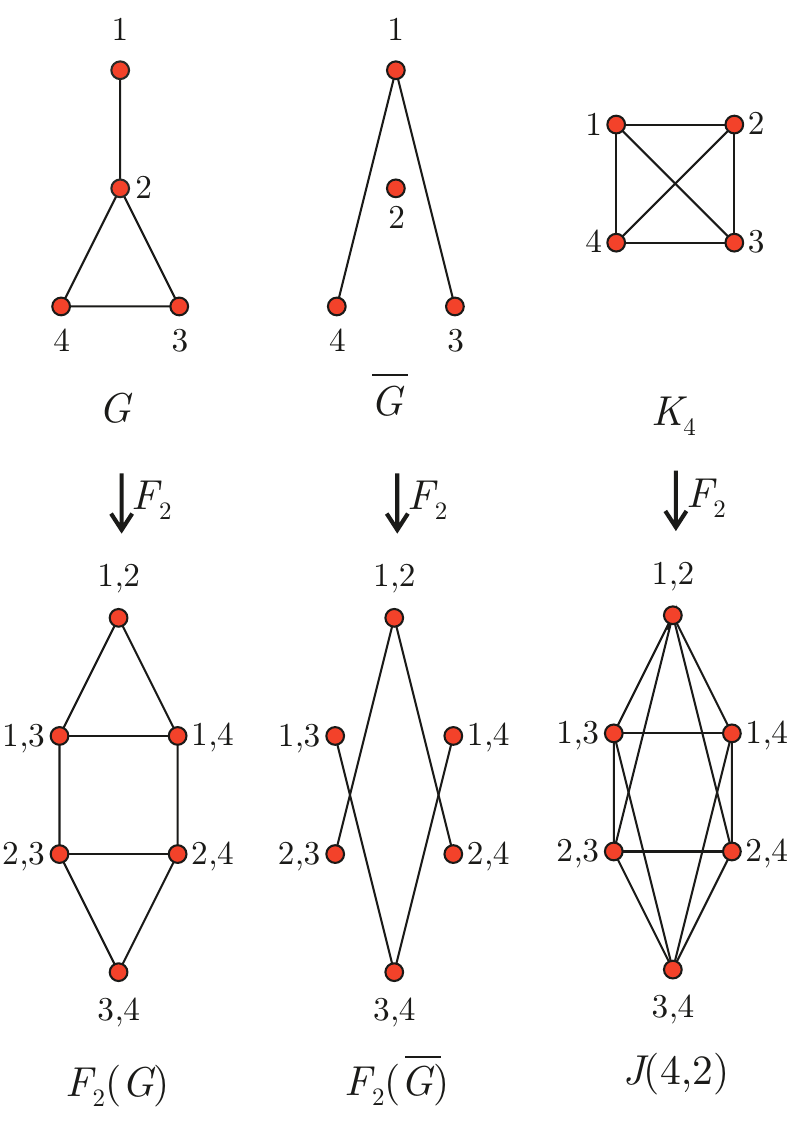}
		\end{center}
		\caption{The graphs $G$, its complement $\overline{G}$, $K_4$, and their 2-token graphs.}
		\label{fig:G+noG+K4}
	\end{figure}

	%%%%%%%%%%%%%%%%%%%%%%%%%%%%%%%%%%%%%%%%%%%%%%%%%%%%%%%%%%%%%%%%%%%%%%
	\begin{example}
		\label{example1}
		In Figure \ref{fig:G+noG+K4}, we show a graph $G$ on 4 vertices, its complement $\overline{G}$, and the complete graph $K_4$. Below them, there are their 2-token graphs.
		Their respective eigenvalues are shown in Table \ref{tab:example1}.
	\end{example}
	
	\begin{table}[!ht]
		\begin{center}
			\begin{tabular}{|c|c|c|c|}
				\hline
				%\noalign{\smallskip}
				Spectrum &  $\spec F_2(G)$ & $\spec F_2(\overline{G})$ & $\spec J(4,2)$\\ 
				\hline\hline
				$\spec F_0 =\spec K_1$	  & 0 & 0 & 0 \\
				\hline
				$\spec F_1\setminus \spec F_0$ & 1 & 3 & 4 \\
				& 3 & 1 & 4 \\
				& 4 & 0 & 4 \\
				\hline
				$\spec F_2\setminus \spec F_1$ & 3  & 3 & 6  \\
				& 5  & 1 & 6  \\
				\hline
			\end{tabular}
		\end{center}
		\caption{The spectra of $F_k(G)$, $F_k(\overline{G})$, and $F_k(K_4)$ for $k=0,1,$ and $2$.}
		\label{tab:example1}
	\end{table}  
	%%%%%%%%%%%%%%%%%%%%%%%%%%%%%%%%%%%%%%%%%%%%%%%%%%%%%%
	
\begin{figure}[!ht]
		\begin{center}
			\includegraphics[width=6cm]{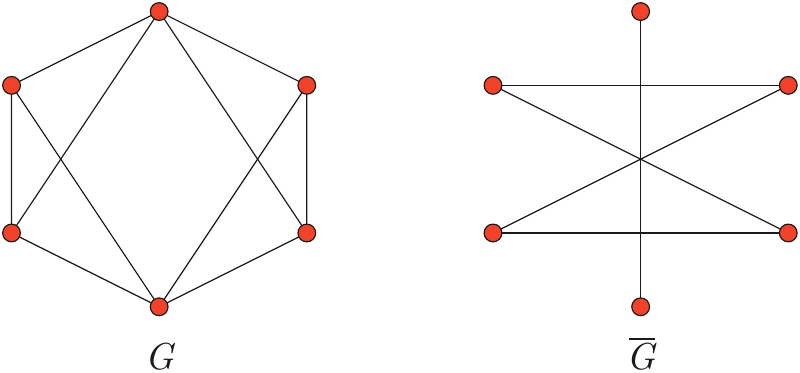}
		\end{center}
		\vskip-0.5cm
		\caption{A graph on 6 vertices and its complement.}
		\label{fig:G+noG+K6}
	\end{figure}

	\begin{example}
		\label{example2}
		Another example is the graph $G$ on 6 vertices shown in Figure \ref{fig:G+noG+K6}, together with its complement.
	\end{example}
	
	\begin{table}[!ht]
		\begin{center}
			%\scriptsize
			%\setlength\tabcolsep{3pt}
			\begin{tabular}{|c|c|c|c|}
				\hline
				Spectrum &  $\spec F_3(G)$ & $\spec F_3(\overline{G})$ & $\spec J(6,3)$ \\
				\hline\hline
				$\spec F_0=\spec K_1$			             & 0 & 0 & 0 \\
				\hline
				& 2 & 4 & 6 \\
				$\spec F_1\setminus \spec F_0$ & 4 & 2 & 6 \\
				& 4 & 2 & 6 \\
				& 4 & 2 & 6 \\
				& 6 & 0 & 6 \\
				\hline
				& 4  & 6 & 10  \\
				& 4  & 6 & 10  \\
				& 6  & 4 & 10  \\
				& 6  & 4 & 10  \\
				$\spec F_2\setminus \spec F_1$ & 6  & 4 & 10 \\
				& 8  & 2 & 10  \\
				& 8  & 2 & 10  \\
				& 8  & 2 & 10 \\
				& 10 & 0 & 10 \\
				\hline
				& 4   & 8  & 12  \\
				& 8   & 4  & 12 \\
				$\spec F_3\setminus \spec F_2$ & 8   & 4  & 12  \\
				& 10  & 2  & 12 \\
				& 10  & 2 & 12 \\
				\hline
			\end{tabular}
		\end{center}
		\caption{The spectra of $F_k(G)$, $F_k(\overline{G})$, and $F_k(K_6)$ for $k=0,1,2,3$.}
	\end{table}
	
	%%%%%%%%%%%%%%%%%%%%%%%%%%%%%%%%%%%%%%%%%%%%%%%%%%%%%%%%%%%%%%
	
	\begin{lemma}
		Let $\B$ be the $(n,k)$-binomial matrix, and let $\A_0(=\I),\A_1,\ldots,\A_k$ be the distance matrices of the Johnson graph $J(n,k)$, for $k\le n-k$.
		Then,
		\begin{equation}
			\M=\B\B^{\top}=\sum_{i=0}^{k-1} (k-i)\A_i.
			%k\I+(k-1)\A_1+(k-2)\A_2+\cdots+\A_{k-1}
			\label{M}
		\end{equation}
	\end{lemma}
	
	\begin{proof}
		We know that each row $\b_r$ of $\B$ represents a vertex $A$ of $J(n,k)=F_k(K_n)$, as the characteristic vector of the 
		$k$-subset $A$. Then, with the row $\b_s$ representing a vertex $B$ such that $\dist(A,B)=i$, we have that, from property {\bf P3}, the vectors $\b_r$ and $\b_s$ have exactly $k-i$ common 1's. Thus,
		$$
		(\M)_{AB}=\sum_{h=1}^n (\B)_{ih}(\B)_{jh}=\b_i\b_j^{\top}=k-i=(k-i)(\A_i)_{AB}.
		$$
		This completes the proof.
	\end{proof}
	
	\begin{corollary}
		\label{coro:commute}
		The matrices $\M$, $\L_k$, $\overline{\L}_k$, and $\L_J$ commute with each other.
	\end{corollary}
	
	\begin{proof}
		Now, we only need to prove that $\M$ commutes with the other matrices. Since $J(n,k)$ is distance-regular, its distance matrices are polynomials of its adjacency matrix $\A$, that is, $\A_i=p_i(\A)$,
		and the same holds for $\M$ because of \eqref{M}. In particular, $J(n,k)$ is  $k(n-k)$-regular and, hence, its Laplacian matrix can be written as 
		\begin{equation}
			\label{LvsA}
			\L_J=\L_k+\overline{\L}_k=k(n-k)\I-\A.
		\end{equation}
		Solving for $\A$, we conclude that $\M$ can be written 
		as a polynomial of $\L_J$, and, therefore, $\M$ commutes with $\L_J$. Finally, since both $\L_k$ and $\overline{\L}_k$ commute with $\L_J$, so does $\M$. 
		% But, for every $i=0,\ldots,k$, the distance graph $G_i$ of $G=J(n,k)$ is $\delta_i$-regular, so that its Laplacian matrices 
	\end{proof}
	
	%%%%%%%%%%%%%%%%%%%%%%%%%%%%%%%%%%%%%%%%%%%%%%%%%%%%%%%%%%55%
	
	\section{A `local' algebra of token graphs}%8
	\label{local-algebra-token}
	
	This section introduces a new algebra generated by the Laplacian matrices $\L_k$ and $\overline{\L}_k$, which contains the Bose-Mesner algebra of Johnson graphs.\\
	Let $\Re[\L_k, \overline{\L}_k]$ be the $\Re$-subalgebra of the $n \times  n$ matrices $M_n(\Re)$ generated by
	the two commuting matrices $\L_k$ and $\overline{\L}_k$. Thus, $\Re[\L_k, \overline{\L}_k]$ consists of all $\Re$-linear combinations of `monomials' $\L_k^i\overline{\L}_k^j$
	where $i$ and $j$ range from
	0 to infinity. Note that $\Re[\L_k, \overline{\L}_k]$ and $M_n(\Re)$ are naturally vector-spaces
	over $\Re$. Moreover, $\Re[\L_k, \overline{\L}_k]$ is a subspace of $M_n(\Re)$.
	
	\begin{theorem}
		\label{th1}
		Let $G$ and $\overline{G}$ be a graph and its complement on $n$ vertices.
		For some $k\le n/2$, let $\L_k$ and $\overline{\L}_k$ be Laplacian matrices  of the token graphs $F_k(G)$ and $F_k(\overline{G})$, respectively.
		Let ${\cal L}(G)=\Re[\L_k, \overline{\L}_k]$ be the $\Re$-vector space of the ${n\choose k} \times {n\choose k}$ matrices $M_n(\Re)$ generated by
		$\L_k$ and $\overline{\L}_k$. Then, the following statements hold:
		\begin{itemize}
			\item[$(i)$]
			${\cal L}(G)$ is a unitary commutative algebra.
			\item[$(ii)$]
			The Bose-Mesner algebra of the Johnson graph $J(n,k)$ is a subalgebra of
			${\cal L}(G)$.
			\item[$(iii)$]
			The dimension of ${\cal L}(G)$ is the number, say $d+1$, of \em{different} pairs $(\lambda_{jr},\overline{\lambda}_{jr})$, for $j=0,\ldots,k$ and $r=1,\ldots,m_j$, defined in Proposition \ref{pro:pairing}.
			\item[$(iv)$]
			If $\dim({\cal L}(G))=d+1$, then there exists a (non-unique) matrix $\R$ such that %$\R=\alpha \L_k+\beta \overline{\L}_k$ such that
			$$
			\{\I,\R,\R^2,\ldots, \R^{d}\}\quad {\mbox and}\quad \{\E_0,\E_1,\E_2,\ldots, \E_{d}\}
			$$
			are bases of ${\cal L}(G)$, where the $\E_i$'s are the idempotents of $\R$.
		\end{itemize}
	\end{theorem}
	
	\begin{proof}
		$(i)$ ${\cal L}(G)$ contains the identity matrix $\I$, and it is commutative since, by Proposition \ref{propo:commute}, 
		$\L_k$ and $\overline{\L}_k$ commute.\\
		%As commented before Proposition \ref{pro:pairing}, the %Laplacian matrix of the $k(n-k)$-regular Johnson graph $J(n,k)$ is $\L_J=\L_k+\overline{\L}_k$.\\
		$(ii)$ From \eqref{LvsA}, the adjacency matrix of the Johnson graphs is $\A_J=n(n-k)\I-\L_J$. Thus,  the Bose-Mesner algebra of $J(n,k)$, generated by $\A_j$, must be a subalgebra of ${\cal L}(G)$.\\
		$(iii)$ This is because, under the hypothesis, we can construct a linear combination of $\L_k$ and $\overline{\L}_k$, say $\R=\alpha\L_k+\beta\overline{\L}_k$, such that the (symmetric) matrix $\R$
		has exactly $d+1$ different eigenvalues. Thus, the minimal polynomial of $\R$ has degree $d+1$, and every matrix of ${\cal L}(G)$ has a minimal polynomial of degree at most $d+1$. Thus, the dimension of 
		${\cal L}(G)$ is $d+1$.\\
		$(iv)$ This follows directly from $(iii)$.
		Notice that, if $\R$ has different eigenvalues $\theta_0,\theta_1,\ldots,\theta_d$, then its idempotents (orthogonal projections onto the $\theta_i$-eigenspaces, see, for instance, Godsil \cite[p. 27]{g93}) can be written as
		$$
		\E_i=\frac{\prod_{j\neq i}(\R-\theta_j\I)}{\prod_{j\neq i}(\theta_i-\theta_j)},\qquad \mbox{for } i=0,1,\ldots,d.
		$$
	\end{proof}
	
	From $(ii)$ and $(iii)$, notice that
	$$
	k+1\le \dim({\cal L}(G))\le {n\choose k}.
	$$
	%\vskip 1cm
	In fact, Gerstenhaber \cite{g61}, as well as Motzkin and Taussky-Todd \cite{mtt95}, proved
	independently that the variety of a commuting pair of matrices $\A$, $\B$ is irreducible so that its dimension is also bounded above by the size of the matrices.\\
	%\vskip 1cm
	The pair $(\L_k,\overline{\L}_k)$ generates an algebraic variety (that is, the collection of all common eigenvectors shared by the two matrices).
	
	\begin{example}
		Considering again Example \ref{example1}, in Table \ref{evR}
		we show the different eigenvalues of the matrix $\R=2\L_{2}+\overline{\L}_2$.
	\end{example}
	
	\begin{table}[H]
		%\begin{flushright}
		\begin{center}
			%\scriptsize
			%\setlength\tabcolsep{3pt}
			\begin{tabular}{|c|c|c|c|}
				\hline
				$\spec \L_2$ & $\spec \overline{\L_2}$ & $\spec \L_J$ & $\spec \R$\\
				\hline\hline
				0 & 0 & 0  & {\bf 0}\\
				\hline
				1 & 3 & 4 & {\bf 5}\\
				3 & 1 & 4 & {\bf 7}\\
				4 & 0 & 4 & {\bf 8}\\
				\hline
				3  & 3 & 6 & {\bf 9} \\
				5  & 1 & 6 & {\bf 11} \\
				\hline
			\end{tabular}
		\end{center}
		%\end{flushright}
		\caption{The different eigenvalues of $\R=2\L_2+\overline{\L}_2$ in Example \ref{example1}.}
		\label{evR}
	\end{table}
	%\vskip-3cm
	%       \begin{figure}[H]
		% 				%\begin{flushleft}
		% \includegraphics[width=4cm]{fig/C_3-amb-aresta+convers+tokens-color.pdf}
		% 					%\end{flushleft}
		%          \end{figure}
	
	Then, for every matrix $\M\in {\cal L}(G)$, there exists a polynomial $p\in \Re^6[x]$ such that $p(\R)=\M$. In particular $H_L(\R)=\J$.
	
\begin{figure}[t]
		\label{L-predist-pols}
		%\vskip-2cm
		\begin{center}
			\includegraphics[width=10cm]{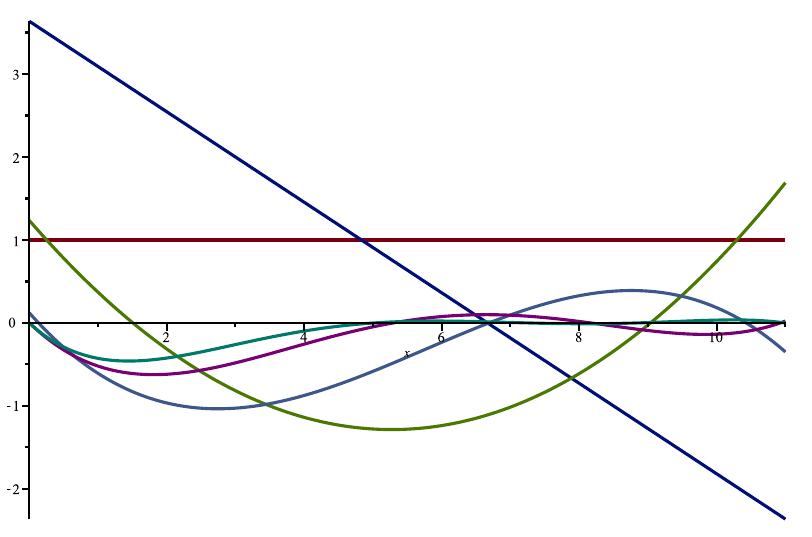}
		\end{center}
		%\vskip-13.5cm
		\caption{The Laplacian polynomials of the matrix $\R$.}
	\end{figure}

	If $\A$ is 1-regular (all eigenvalues are different), then any matrix $\B$ that commutes
	with $\A$  must be a polynomial in $\A$. Described differently, $\B$ is already in
	the algebra $\Re[\A]$, that is, $\Re[\A, \B] = \Re[\A]$. However, $\Re[\A]$ is of dimension $n$ as $\A$
	is 1-regular, so $\Re[\A,\B]$ has dimension $n$.
	The Laplacian predistance polynomials of $\R$ are (see Figure \ref{L-predist-pols} for a plot of them):
	\begin{align*}
		q_0(x) &=1,\\
		q_1(x) &=\textstyle -\frac{6}{11}x+\frac{40}{11},\\
		q_2(x) &=\textstyle \frac{75}{827}x^2 - \frac{11250}{9097}x  +\frac{8700}{9097},\\
		q_3(x) &=\textstyle -\frac{3255}{248788}x^3 + \frac{283185}{1254559}x^2 - \frac{193990839}{205747676}x + \frac{6357015}{51436919},\\
		q_4(x) &=\textstyle \frac{424}{248327}x^4-\frac{643258880}{15445194419}x^3+\frac{122491480}{376712059}x^2-\frac{12482051200}{15445194419}x+\frac{51775488}{15445194419},\\
		q_5(x) &=\textstyle -\frac{1}{4620}x^5+\frac{398710}{57363537}x^4-\frac{891947}{10926388}x^3+\frac{23788868}{57363537}x^2-\frac{439657769}{573635370}x+\frac{6}{248327}. \\
	\end{align*}
	
	\section{A `global' algebra of token graphs}
	\label{global-algebra-token}
	
	For given values of $n$ and $k$, the `local' algebra of Section \ref{local-algebra-token} is constructed from the Laplacian matrices of the $k$-token graphs of $G$ and its complement $\overline{G}$. In this section, we only fixed $n$ and present an algebra containing the Laplacian and adjacency matrices of {\bf all} the $k$-token graphs of a graph on $n$ vertices.
	
	Let $G$ be a graph with $n$ vertices and $m$ edges. Notice that every edge $e=\{u,v\}$ of $G$ gives rise to ${n-2\choose k-1}$ edges
	in the $k$-token graph $F_k(G)$. Indeed, we have an edge in $F_k(G)$ every time there is one token moving between $u$ and $v$, and the other $k-1$ tokens are in some fixed position in $V\setminus \{u,v\}$. In particular, each of the ${n\choose 2}$ edges of $G=K_n$, produces ${n-2\choose k-1}$ edges in the Johnson graph $J(n,k)=F_k(K_n)$.  This suggests the following edge-decomposition of $J(n,k)$. 
	For every edge $e$ of $K_n$, we consider the graph $K_n(e)$ that consists of the single edge $e$ plus the remaining $n-2$ vertices of $K_n$ with no edges between them. Then, we get the following result.
	
	\begin{lemma}
		\label{lem:Le}
		\begin{itemize}
			\item[$(i)$] For every edge $e$ of $K_n$, the $k$-token graph $F_k(K_n(e))$ consists of 
			${n-2\choose k-1}$ independent edges plus ${n\choose k}-2{n-2\choose k-1}$ isolated vertices.
			\item[$(ii)$] The edges of the Johnson graph $J(n,k)$ are the union of the edges
			of $F_k(K_n(e))$ when $e$  runs through all the edges of $K_n$. So,
			the edges of the $k$-token graph $F_k(G)$ of a graph $G\subset K_n$ are the union of the edges
			of $F_k(K_n(e))$ when $e$  runs through all the edges of $G$.
			\item[$(iii)$] The Laplacian [adjacency] matrix $\L_J$ [$\A_J$] of $J(n,k)$ is the sum of the Laplacian [adjacency] matrices $\L_e$ [$\A_e$] of
			$F_k(K_n(e))$. 
			\item[$(iv)$] 
			Let $\L_G$ [$\A_G$] be the Laplacian [adjacency] matrix of a spanning subgraph $G(n,k)$ of $J(n,k)$. Then $G(n,k)$ is the $k$-token graph of some graph $H$ on $n$ vertices if and only if $\L_G$ [$\A_G$]  is the sum of the Laplacian [adjacency] matrices $\L_e$ [$\A_e$] of
			$F_k(K_n(e))$ where $e$ runs through all edges of $H$:
			\begin{equation}
				\label{eq:LGvsLe}
				\L_G=\sum_{e\in E(H)} \L_e,\qquad
				\A_G=\sum_{e\in E(H)} \A_e.
			\end{equation}
		\end{itemize}
	\end{lemma}
	
	\begin{proof}
		$(i)$ Here we only need to prove that all the edges of $F_k(K_n(e))$ have no vertices in common. If $e=\{a,b\}$ is an edge of $K_n$, then two any edges of $F_k(K_n(e))$ are of the form $\{A,B\}$, $\{C,D\}$, with $A=A'\cup \{a\}$, $B=A'\cup \{b\}$, $C=A''\cup \{a\}$, and
		$D=A''\cup \{b\}$, where $|A'|=|A''|=k-1$ and $A'\neq A''$. Thus, $A,B,C,D$ are different vertices of $F_k(K_n(e)$. $(ii)$--$(iv)$ All the statements in these items follow from the following fact: Let $G$ be a subgraph of $K_n$. Then, there exist an edge $\{A,B\}$ of $F_k(G)$ if and only if $A\Delta B=\{a,b\}$ for some edge $e=\{a,b\}$ in $G$ (and also in $K_n$).
		Thus, for $A\neq B$, $(\L_G)=-1$ if and only if $(\L_e)_AB=-1$. Moreover, if $A=B$, all the edges $e$ of $F_k(G)$ with endpoint $A$ contribute with the value $(\L_e)=1$ and, hence $(\L_G)_AA=\deg_{F_k(G)} A$. The equality for the adjacency matrix $\A_G$ is also clear.
	\end{proof}
	
	In other words, the $k$-token graph of a graph $H$ on $n$ vertices can be seen as the sum of the $k$-token graphs of some `elementary' graphs $K_n(e)$.  In terms of matrices, the Laplacian or adjacency matrices of $F_k(G)$ are a linear combination, with coefficients $0$ and $1$ of the `elementary' matrices $\L_e$ or $\A_e$, respectively.
	
	In particular, if a Laplacian matrix can be written as a linear combination of the matrices $\L_e$, $e\in K_n(e)$, then $\L_e$ commutes with $\L_J$ (since $\L_e$ is, in fact the Laplacian matrix of a $k$-token graph, and Corollary \ref{coro:commute} applies).
	
	\begin{theorem}
		For some fixed $n$ and $k(\le n/2)$, let $\A_e$ be, for $e\in E(K_n)$, the adjacency matrices  of the token graphs $F_k(K_n(e))$.
		Let ${\cal A}(n,k)$ be the  $\Re$-vector space of the ${n\choose k} \times {n\choose k}$ matrices $M_n(\Re)$ generated by
		the matrices $\A_e$. Then, the following statements hold:
		\begin{itemize}
			\item[$(i)$]
			${\cal A}(n,k)$ is a unitary non-commutative algebra.
			\item[$(ii)$]
			For any graph $G\subset K_n$, the algebra ${\cal L}(G)$ of Theorem \ref{th1} is a subalgebra of
			${\cal A}(n,k)$.
			\item[$(iii)$]
			The dimension of ${\cal A}(n,k)$ is ${n\choose 2}$, with a basis constituted by the matrices $\A_e$, for $e\in E(K_n)$.
		\end{itemize}
	\end{theorem}
	
	\begin{proof}
		$(i)$ It is clear since, if the edges $e$ and $e'$ has a vertex in common,  the matrices $\A_e$ and $\A_{e'}$ do not commute.  $(ii)$ Using \eqref{eq:LGvsLe}, we only need to prove that the Laplacian matrices $\L_e$, for $e\in K_n(e)$, belong to  ${\cal A}(n,k)$. As $\A_e$ corresponds to a graph with a series of independent vertices, $\A_e^2$ is a diagonal matrix with 1's in the entries corresponding to the endpoints of such edges. Then, the statement follows from the equality $\L_e=\A_e^2-\A_e$ (see Example \ref{Ae-Le}).
		$(iii)$ No pair of matrices $\A_e$ have a common entry $1$, so that they are linearly independent.
	\end{proof}
	
	\begin{example}
		\label{Ae-Le}
		In the case of $K_4$, with vertices $1,\ldots,4$, the 
		% ${4\choose 2}\times {4\choose 2}$ 
		adjacency and Laplacian matrices, indexed by the vertices of $F_2(K_4)$ in the order $12,13,14,23,24,34$, and induced by the edge $e=\{1,2\}$  are:
		
		$$
		\A_{e} = 
		\left(
		\begin{array}{cccccc}
			0 & 0 & 0 & 0 & 0 & 0\\
			0 & 0 & 0 & 1 & 0 & 0\\
			0 & 0 & 0 & 0 & 1 & 0 \\
			0 & 1 & 0 & 0 & 0 & 0\\
			0 & 0 & 1 & 0 & 0 & 0\\
			0 & 0 & 0 & 0 & 0 & 0
		\end{array}
		\right), \quad
		% \A_{13} = 
		% \left(
		% \begin{array}{cccccc}
			% 0 & 0 & 0 & 1 & 0 & 0\\
			% 0 & 0 & 0 & 0 & 0 & 0\\
			% 0 & 0 & 0 & 0 & 0 & 1 \\
			% 1 & 0 & 0 & 0 & 0 & 0\\
			% 0 & 0 & 0 & 0 & 0 & 0\\
			% 0 & 0 & 1 & 0 & 0 & 0
			% \end{array}
		% \right)
		% $$
		% $$
		\L_{e} = \A_{e}^2-\A_{e}=
		\left(
		\begin{array}{cccccc}
			0 & 0 & 0 & 0 & 0 & 0\\
			0 & 1 & 0 & -1 & 0 & 0\\
			0 & 0 & 1 & 0 & -1 & 0 \\
			0 & -1 & 0 & 1 & 0 & 0\\
			0 & 0 & -1 & 0 & 1 & 0\\
			0 & 0 & 0 & 0 & 0 & 0
		\end{array}
		\right). 
		% \L_{13} = 
		% \left(
		% \begin{array}{cccccc}
			% 1 & 0 & 0 & -1 & 0 & 0\\
			% 0 & 0 & 0 & 0 & 0 & 0\\
			% 0 & 0 & 1 & 0 & 0 & -1 \\
			% -1 & 0 & 0 & 1 & 0 & 0\\
			% 0 & 0 & 0 & 0 & 0 & 0\\
			% 0 & 0 & -1 & 0 & 0 & 1
			% \end{array}
		% \right)
		$$
	\end{example}
	
	Concerning the Laplacian matrices, the following result shows that the algebra 
	${\cal L}(G)$ of Theorem \ref{th1} can be constructed only if the condition of commutativity in Proposition \ref{propo:commute} holds.
	\begin{proposition}
		Let $\L_k$ and $\overline{\L}_k$ be, respectively, the Laplacian matrices of a spanning subgraph $G(n,k)$ of $J(n,k)$ and its complement $\overline{G}(n,k)$ with respect to $J(n,k)$. Then, $\L_k$ and $\overline{\L}_k$ commute if and only if  $G(n,k)$ and  $\overline{G}(n,k)$ are the $k$-token graphs of some graph $H$, on $n$ vertices, and its complement $\overline{H}$.
	\end{proposition}
	
	\begin{proof}
		If $G(n,k)=F_k(H)$ and $\overline{G}(n,k)=F_k(\overline{H})$, then we know, by Proposition \ref{propo:commute}, that $\L_k$ and $\overline{\L}_k$ commute.
		To prove the converse, we use Corollary \ref{coro:commute}, by noting that $\L_k$ and $\overline{\L}_k$ commute if and only if $\L_k$ and $\L_J$ do.
		However, if $\L_k$ and $\L_J$ do not commute,
		$\L_k$ cannot be written as a $(0,1)$-linear combination of some matrices $\L_e$, for $e\in K_n(e)$ (see the last comment after Lemma \ref{lem:Le}).
		Consequently, by the same lemma, 
		if $\L_k$ and $\overline{\L}_k$ do not commute, neither $G(n,k)$ nor $\overline{G}(n,k)$ can be $k$-token graphs.
	\end{proof}
	% \section{Problems and (possible) future work}
	% \begin{itemize}
		% \item
		% Use this algebra in the context of codes or designs.
		% \item
		% Prove (or disprove) that all distance-regular graphs with the same parameters have cospectral 2-token (symmetric square) graphs.
		% \item
		% What about $\L_1+\L_2+\L_3=\L_J$?
		% \item
		% Consider the case when $G$ and $\overline{G}$ are strongly regular or self-complementary graphs, as the Payley graphs.
		% \end{itemize}
	
	% %\textbf{Payley graphs}
	% \begin{figure}[H]
		% 	\begin{center}
			% \includegraphics[width=7cm]{fig/Payley.jpg}
			% 	\end{center}
		% \end{figure}


\begin{thebibliography}{10}
		\label{bibliography}
		\bibitem{abel91}
		Y. Alavi, M. Behzad, P. Erd\H{o}s, and D. R. Lick,
		Double vertex graphs,
		{\em J. Comb. Inf. Syst. Sci.} {\bf 16} (1991), no. 1, 37--50.
		
		\bibitem{all02}
		Y. Alavi, D. R. Lick, and J. Liu, 
		Survey of double vertex graphs, 
		\emph{Graphs Combin.}  \textbf{18} (2002) 709--715.
		
		\bibitem{acss12}
		B. Arsi\'{c}, D. Cvetkovi\'{c}, S. K. Simi\'{c}, and M. \v{S}kari\'{c},
		Graph spectral techniques in computer sciences, 
		\textit{Appl. Anal. Discrete Math.} \textbf{6} (2012), no. 1, 1--30.
		
		\bibitem{agrr07}
		K. Audenaert, C. Godsil, G. Royle, and T. Rudolph,
		Symmetric squares of graphs,
		\emph{J. Combin. Theory B} \textbf{97} (2007) 74--90.
		
		%\bibitem{bns78}
		%J. Bunch, C. Nielsen, and D. Sorensen,
		%Rank-one modificaton of the symmetric eigenproblem,
		%\emph{Numer. Math.} \textbf{31} (1978) 31--48.
		
		%\bibitem{clr10}
		%P. Caputo, T. M. Liggett, and T. Richthammer, 
		%Proof of Aldous' spectral gap conjecture,
		%{\em J. Amer. Soc.} {\bf 23} (2010), no. 3, 831--851.
		
		\bibitem{cflr17}
		W. Carballosa, R. Fabila-Monroy, J. Lea\~nos, and L. M. Rivera, 
		Regularity and planarity of token graphs, 
		{\em Discuss. Math. Graph Theory} {\bf 37} (2017), no. 3, 573--586.
		
		%\bibitem{c16}
		%F. Cesi,
		%A few remarks on the octopus inequality and Aldous' spectral gap conjecture,
		%{\em Comm. Algebra} {\bf 44} (2016), no. 1, 279--302.
		
		% \bibitem{cr12}
		% E. S. Coakley and V. Rokhlin, 
		% A fast divide-and-conquer algorithm for computing the spectra of real symmetric tridiagonal matrices, 
		% \textit{Appl. Comput. Harmon. Anal.} \textbf{34} (2013) 379--414.
		
		\bibitem{ddffhtz21}
		C. Dalf\'o, F. Duque, R. Fabila-Monroy, M. A. Fiol, C. Huemer, A. L. Trujillo-Negrete, and F. J. Zaragoza Mart\'{\i}nez, 
		On the Laplacian spectra of token graphs, 
		{\em Linear Algebra Appl.} {\bf 625} (2021) 322--348.
		
		\bibitem{df22}
		C. Dalf\'o and M. A. Fiol, 
		On the algebraic connectivity of token graphs, 
		\textit{J. Algebraic Combin.}, accepted, 2023, \texttt{https://arxiv.org/abs/2209.01030}.
		
		%\bibitem{dfm22}
		%C. Dalf\'o, M. A. Fiol, and A. Messegu\'e, Some bounds on the algebraic connectivity of token graphs, submitted (2022).
		
		%\bibitem{dfmrs17}
		%C. Dalf\'o, M. A. Fiol, M. Miller, J. Ryan, and J. \v{S}ir\'a\v{n},
		%An algebraic approach to lifts of digraphs,
		%{\em Discrete Appl. Math.} {\bf 269} (2019) 68--76.
		
		%\bibitem{df22b}
		%C. Dalf\'o, M. A. Fiol, S. Pavl\'ikov\'a, and J. \v{S}ir\'an, 
		%On the spectra and eigenspaces of the universal adjacency matrices of arbitrary lifts of graphs, 
		%\textit{Linear Multilinear Algebra} {\bf 71} (2023), no. 5, 693--710.
		
		%\bibitem{dfs19}
		%C. Dalf\'o, M. A. Fiol, and J. \v{S}ir\'a\v{n},
		%The spectra of lifted digraphs, {\em J. Algebraic Combin.} \textbf{50} (2019) 419--426.
		
		\bibitem{vf14}
		E. R. van Dam and  M. A. Fiol,
		The Laplacian spectral excess theorem for distance-regular graphs,
		{\em Linear Algebra Appl.} {\bf 458} (2014) 245--250.
		
		%\bibitem{f23}
		%R. Fabila-Monroy, Personal communication.
		
		\bibitem{ffhhuw12}
		R. Fabila-Monroy, D. Flores-Pe\~{n}aloza, C. Huemer, F. Hurtado, J. Urrutia, and D. R. Wood,
		Token graphs,
		\emph{Graphs Combin.} \textbf{28} (2012), no. 3, 365--380.
		
		%\bibitem{fi73}
		%M. Fiedler,
		%Algebraic connectivity of graphs,
		%{\em Czech. Math. Journal} {\bf 23} (1973), no. 2, 298--305.
		
		\bibitem{fg97}
		M. A. Fiol and E. Garriga, 
		From local adjacency polynomials to locally pseudo distance-regular graphs, 
		{\em J. Combin. Theory Ser. B} {\bf 7I} (1997) 162--183.
		
		\bibitem{fgy96}
		M. A. Fiol, E. Garriga, and J. L. A. Yebra, 
		Locally pseudo-distance-regular graphs, 
		{\em J. Combin. Theory Ser. B} {\bf 68} (1996) 179--205.
		
		% \bibitem{f}
		% G. Frobenius,
		% Ueber lineare Substitutionen und bilineare Formen, {\em Journal f\"{u}r die reine und angewandte Mathematik} {\bf 84} (1877) 1--63.
		
		\bibitem{g61}
		M. Gerstenhaber, 
		On dominance and varieties of commuting matrices, 
		{\em Ann. of Math.} {\bf 73} (1961) 324--348.
		
		\bibitem{g93}
		C. D. Godsil,
		Algebraic Combinatorics,
		Chapman and Hall, New York, 1993.
		
		% \bibitem{g95}
		% C. D. Godsil,
		% Tools from linear algebra, in {\em Handbook of Combinatorics} (eds. Graham, Gr\"otschel, Lov\'asz), MIT press 1995, pp. 1705--1748.
		
		%\bibitem{gr01}
		%C. Godsil and  G. Royle,
		%{\em Algebraic Graph Theory},
		%Graduate Texts in Mathematics {\bf 207}, Springer-Verlag, New York, 2001.
		
		%\bibitem{go99}
		%H. W. Gould, The Girard-Waring power sum formulas for symmetric functions and Fibonacci sequences, {\em Fibonacci Quart.} {\bf 37} (1999), no. 2, 135--140.
		
		% \bibitem{gm94}
		% R. Grone and R. Merris, The Laplacian spectrum of a graph II, {\em SIAM J. Discrete Math.} {\bf 7}
		% (1994), 221--229.
		
		%\bibitem{gms90}
		%R. Grone, R. Merris, and V. S. Sunder, The Laplacian spectrum of a graph, {\em SIAM J. Matrix Anal. Appl.} {\bf 11} (1990) 218–238.
		
		% \bibitem{gt77}
		% J. L. Gross and T. W. Tucker,
		% Generating all graph coverings by permutation voltage assignments,
		% {\em Discrete Math.} {\bf 18} (1977) 273--283.
		
		\bibitem{h63}
		A. J. Hoffman, 
		On the polynomial of a graph, 
		{\em Amer. Math. Monthly} {\bf 70} (1963) 30--36.
		
		%\bibitem{ir22}
		%S. Ibarra and L. M. Rivera,
		%The automorphism groups of some token graphs,
		%{\tt arXiv:1907.06008v3[math.CO]}
		
		%\bibitem{k10}
		%S. Kirkland,
		%Algebraic connectivity for vertex-deleted subgraphs, and a notion of vertex centrality,
		%{\em Discrete Math.} {\bf 310} (2010) 911--921.
		
		\bibitem{mtt95}
		T. Motzkin and O. Taussky-Todd, 
		Pairs of matrices with property L. II, 
		{\em Trans. Amer. Math. Soc.} {\bf 80} (1955) 387--401.
		
		\bibitem{nk07}
		V. Nikiforov, The spectral radius of subgraphs of regular graphs,
		{\em Electron. J. Combin.} {\bf 14} (2007) \#N20.
		
		%\bibitem{l23}
		%A. Lew,
		%Garland’s method for token graphs,
		%{\tt https://arxiv.org/abs/2305} {\tt .02406v1}, 2023.
		
		%\bibitem{o19}
		%Y. Ouyang, Computing spectral bounds of the Heisenberg ferromagnet from geometric considerations, {\em J.  Math. Physics} {\bf 60} (2019) 071901.
		
		% \bibitem{pl08}
		% K. L. Patra and A. K. Lal,
		% The effect on the algebraic connectivity of a tree by grafting or collapsing of edges,
		% {\em Linear Algebra Appl.} {\bf 428} (2008) 855--864.
		
		%\bibitem{ps17}
		%K. L. Patra and B. K. Sahoo,
		%Bounds for the Laplacian spectral radius of graphs,
		%{\em Electron. J. Graph Theory and Appl.} {\bf 5} (2017), no. 2, 276--303.
		
		%\bibitem{y05}
		%W.-C. Yueh, 
		%Eigenvalues of several tridiagonal matrices, 
		%\textit{Appl. Math. E-Notes} \textbf{5} (2005) 66--74.
	\end{thebibliography}
\end{document}